\documentclass[a4paper,11pt,leqno]{amsart}
\usepackage[big]{opposite_filtrations}
\usepackage{xcolor,fancyhdr}

\newcommand{\C}{\mathbb C}\newcommand{\J}{\mathcal J}
\newcommand{\tJ}{\widetilde{\mathcal J}}

\newcommand{\Span}{\operatorname{span}}

\newcommand{\eps}{\epsilon}\newcommand{\dd}{\delta}

\title{Opposite filtrations on the parabolic compactified Jacobian}
\author{Oscar Kivinen}
\address{Aalto University, Department of Mathematics and Systems Analysis, P.O. Box 11100, FI-00076 Aalto, Finland}
\email{oscar.kivinen@aalto.fi}
\begin{document}
\maketitle
\begin{abstract}
We show that the Hikita filtration on the cohomology of the parabolic compactified Jacobian of the plane curve singularity $x^m=y^n$, $(m,n)=1$, is opposite to the perverse filtration in a precise sense. 
This implies the parabolic Hilb-vs-Quot conjecture from \cite{KTpublished}, and consequently the main conjecture from \cite{ORS} for these singularities.
The new ingredient is an analysis of how the {\em Bott--Mitchell filtration} on the homology of the affine Grassmannian of $SL_n$ relates to both of the above filtrations.
\end{abstract}
\tableofcontents

\section{Introduction}
Over the last 15 years, the ORS conjecture \cite[Conjecture 2]{ORS} has spurred a lot of activity around the algebraic geometry of various moduli spaces of sheaves on singular curves, as well as the topological invariants of the corresponding links of the singularities. The simplest singularities are the quasihomogeneous plane curve singularities $C=\{y^n=x^m\}\subset \bb{A}^2_\bb{C}$, $(n,m)=1$, whose links are exactly the positive torus knots. 
As explained in \cite[\S1.6]{KTpublished}, the ORS conjecture was verified for the $n=2$ series of these singularities in the original paper. In \cite[Theorem 5]{KTpublished}, we proved a slightly modified version of the original conjecture, by replacing the Hilbert schemes by Quot-schemes of the pushforward of the structure sheaf from the normalization. 
This also allowed us to verify the original conjecture for the $n=3$ series of quasihomogeneous singularities. The conjectural general equivalence of the Quot scheme and Hilbert scheme computations is the eponymous Hilb-vs-Quot conjecture.

In this note, we prove a certain stronger statement about two filtrations on the cohomology of the parabolic compactified Jacobian for the quasihomogeneous germs (Theorem \ref{thm:coprime}). As explained in the first version of \cite{KTpublished} and recalled here, this implies the parabolic Hilb-vs-Quot conjecture \cite[Conjecture 2]{KTpublished} and in combination with \cite[Theorem 5(i)]{KTpublished}, also implies \cite[Conjecture 2]{ORS} in this case. We also note that our results imply \cite[Conjecture A.12]{ORS}. The filtrations in question are the {\em gap filtration} from \cite[\S8.2]{KTpublished}, which was explained to coincide up to an involution with the {\em Hikita filtration}, as well as the {\em perverse filtration}.

We follow the notation of \cite[\S1 and \S8]{KTpublished}, with
$f:=y^n-x^m$. If $\nu$ is a composition of $n$, write
$\J^\nu:=\CptJac_\nu\cong\CptPic_\nu/\Gamma$ for the parabolic
compactified Jacobian of $R=\C[[x,y]]/(f)$, and write
$\J:=\J^{(n)}$ and $\tJ:=\J^{(1,\ldots,1)}$ for the unflagged and
full-flag cases. All homology and cohomology groups have coefficients in a field containing $\bb{Q}$, we take complex coefficients for simplicity. The closed subvariety $\J^\nu_{\le r}$ consists of points
whose gap invariant is at most $r$. We define
\[
 \mathsf{Q}^{>r}H^*(\J^\nu)
 :=\ker\bigl(H^*(\J^\nu)\longrightarrow H^*(\J^\nu_{\le r})\bigr), \;\; \gr^{\mathsf{Q}}_rH^*(\J^\nu):=\mathsf{Q}^{>r-1}H^*(\J^\nu)/\mathsf{Q}^{>r}H^*(\J^\nu).
\]
The increasing perverse filtration $\mathsf{P}_{\le p}$ is the one obtained
from the parabolic Hitchin fibration, as in
\cite[Definition 8.5.1]{OY16}; in the unflagged case see also
\cite[\S5.3]{OY17}. We use the normalization $1\in \mathsf{P}_{\le0}$, in
which the grading on the associated Cherednik module is polynomial
degree. For partial flags this filtration agrees, by Springer theory, with the $S_\nu:=S_{\nu_1}\times\cdots\times S_{\nu_{\ell(\nu)}}$-invariant part of the full-flag
filtration.

With this notation in place, our main theorem is the following, which is a parabolic analog of a stronger version of \cite[Conjecture 8.10]{KTpublished}.\footnote{The inspiration for this result is in fact Conjecture 12 from the first arXiv version of \cite{KTpublished}, which was cut from the published version.}
\begin{thm}\label{thm:coprime}
Let $\nu$ be a composition of $n$ and $j,k$ two integers. Then, the $\mathsf{Q}$- and $\mathsf{P}$-filtrations oppose each other in the following sense:
\begin{equation}
\label{eq:coprime}
 H^{2k}(\J^\nu)=\mathsf{P}_{\le k+j}H^{2k}(\J^\nu)
                  \oplus \mathsf{Q}^{>j}H^{2k}(\J^\nu),
 \qquad H^{\mathrm{odd}}(\J^\nu)=0.
\end{equation}
In particular, the two filtrations give canonical isomorphisms
\[
 \gr^{\mathsf{P}}_{k+j}H^{2k}(\J^\nu)
 \cong \mathsf{P}_{\le k+j}H^{2k}(\J^\nu)\cap \mathsf{Q}^{>j-1}H^{2k}(\J^\nu)
 \cong \gr^{\mathsf{Q}}_jH^{2k}(\J^\nu).
\]
\end{thm}

\subsection*{Acknowledgements and AI disclosure}
Theorem \ref{thm:coprime} was obtained in discussion with GPT-6 Astra, which also helped us with checking examples and writing.
The argument presented here is a simplified and conceptually clarified version of the original one.
We thank Minh-Tam Trinh for comments on a draft of this paper.
\section{The proof of the main theorem}
\subsection{The Bott--Mitchell filtration}
We follow the notation of \cite[\S1 and \S7--\S8]{KTpublished}, with
$f:=y^n-x^m$. In particular, we let $S=\C[[\varpi]]$, $R=\C[[\varpi^n, \varpi^m]]$ and denote $x=\varpi^n$, and $y=\varpi^m$.
The delta invariant is $\dd=\dim_\C(S/R)=(m-1)(n-1)/2$.

We identify the affine Grassmannian $\Grass_{GL_n}$ with the space of $\C[[x]]$-lattices in $\C((\varpi))$ by fixing the $\C[[x]]$-basis $1,\varpi,\varpi^2,\ldots,\varpi^{n-1}$. 
Let $\mathcal G$ be the component of $\Grass_{GL_n}$ consisting of lattices relative index $\delta$ with respect to $S$, 
identified with $\Grass_{SL_n}$
by multiplication by $\varpi^{-\dd}$. We realize $\J$ as the closed
subvariety of $y$-stable $\C[[x]]$-lattices in $\mathcal G$.
Every such lattice $M$ lies in $S$ and has colength $\dd$ there.

For $r\ge0$, define the projective Schubert variety
\begin{equation}\label{eq:gap-quot}
 Z_r:=\{M\in\mathcal G:M\subseteq\varpi^{\dd-r}S\}
 \xrightarrow[\sim]{\ M\mapsto\varpi^{r-\dd}M\ }
 \operatorname{Quot}^{r}_{\C[[x]]}(S).
\end{equation}
The Quot scheme on the right parametrizes colength-$r$
$\C[[x]]$-submodules of $S$; as above, it is identified with the $GL_n$ Schubert variety of
coweight $(r,0,\ldots,0)$. The varieties $Z_r$ form an increasing
exhaustion of $\mathcal G$. Since the gap invariant is
$c(M)=\dd-\min\operatorname{val}(M)$ in this realization,
\[
 \J_{\le r}=\J\cap Z_r
 \cong\operatorname{Quot}^{r}_R(S).
\]
Thus the gap filtration is induced by these Schubert varieties; see
\cite[Lemma 8.1 and Theorem 8.3]{KTpublished} for its relation to Hikita's filtration.

We write $H_*(\mathcal G)$ for the direct limit of the homology of
its finite-dimensional projective Schubert varieties; ordinary and
Borel--Moore homology agree on each of them. With its convolution product,
\[
 H_*(\mathcal G)\cong \C[u_1,\ldots,u_{n-1}],\qquad \deg u_i=2i.
\]
Here $u_i$ is the image of the class of a linear $\mathbb P^i$ in
$Z_1\cong\mathbb P^{n-1}$. Define the increasing filtration
\[
 \mathsf{F}_{\le r}H_*(\mathcal G)
 :=\operatorname{im}\bigl(H_*(Z_r)\longrightarrow H_*(\mathcal G)\bigr).
\]
Mitchell's filtration theorem \cite{Mitchell}, recalled in
\cite[\S1 and \S11]{BM}, states that these homology maps are injective
and that this is the filtration by the naive polynomial degree in the $u_i$:
\begin{equation}\label{eq:homology-length}
 \mathsf{F}_{\le r}H_*(\mathcal G)
 =\Span\{u_1^{\alpha_1}\cdots u_{n-1}^{\alpha_{n-1}}:
                 \alpha_1+\cdots+\alpha_{n-1}\le r\}.
\end{equation}
Note that the
homological degree of the monomials above is $2\sum_i i\alpha_i$. The dual decreasing
filtration is
$\mathsf{F}^{>r}H^*(\mathcal G)=\ker(H^*(\mathcal G)\to H^*(Z_r))$.

Let $u_i^\vee\in H^{2i}(\mathcal G)$ be dual to the monomial $u_i$
in the monomial basis of $H_*(\mathcal G)$: it pairs to $1$ with $u_i$
and to $0$ with every other monomial of degree $2i$. These classes
generate $H^*(\mathcal G)=\C[u_1^\vee,\ldots,u_{n-1}^\vee]$.
Here and below the pairing is ordinary homology-cohomology duality.

\begin{lem}\label{lem:bott}
Restriction induces an isomorphism of vector spaces
\[
 \Span\left\{\prod_{i=1}^{n-1}(u_i^\vee)^{\alpha_i}:
       \sum_i i\alpha_i=k,\quad\sum_i\alpha_i\le r\right\}
 \xrightarrow{\sim}H^{2k}(Z_r).
\]
On the Quot scheme in \eqref{eq:gap-quot}, let $V$ be the rank-$r$
universal quotient bundle, with fiber $S/N$ at a lattice $N$.
Then $u_i^\vee|_{Z_r}=p_i(V)$, where $p_i(V)$ is the $i$th power
sum of the Chern roots of $V$.
\end{lem}
\begin{proof}
The quotient bundles fit into exact sequences under convolution, so
their Chern power sums are additive. On $Z_1=\mathbb P^{n-1}$ the
quotient bundle is $\mathcal O(1)$. Thus $p_i(V)$ pairs to $1$ with
$u_i$ and vanishes on products of positive-degree homology classes;
this identifies it with the specified dual class.

For multi-indices $\alpha,\beta$, write $|\alpha|=\sum_i\alpha_i$
and $u^\beta=\prod_i u_i^{\beta_i}$. Compatibility of convolution
with the pairing gives
\[
 \left\langle\prod_i(u_i^\vee)^{\alpha_i},u^\beta\right\rangle=0
       \quad\text{if }|\beta|>|\alpha|,
\]
and, at equal lengths, this pairing is zero unless $\alpha=\beta$,
when it equals $\prod_i\alpha_i!$. In each degree, the monomials
of length at most $r$ on the two sides therefore have a square
triangular pairing matrix with nonzero diagonal. Now use
\eqref{eq:homology-length}.
\end{proof}
\subsection{Chern classes and restriction}
Let $E$ be the rank-$n$ vector bundle on $\mathcal G$ with fiber
$M/xM$ over $M$. Its full flag bundle $\operatorname{Fl}(E)$ is the component of the affine flag variety above $\mathcal G$. 
A point 
is a chain
$M=M_0\supset M_1\supset\cdots\supset M_n=xM$ with
$\dim(M_{j-1}/M_j)=1$. Denote the quotient line bundle with this
fiber by $L_j$, and put $\xi_j=c_1(L_j)$. The subvariety $\tJ$
consists of the $y$-stable chains, and $\tJ_{\le r}$ is its
intersection with $\operatorname{Fl}(E|_{Z_r})$.

On the Quot scheme in \eqref{eq:gap-quot}, we use $N/xN$ and $S/xN$
to denote the vector bundles with these fibers at a lattice $N$;
$S/xS$ denotes the constant rank-$n$ bundle. There are exact
sequences of ordinary vector bundles
\[
 \begin{gathered}
 0\longrightarrow N/xN\longrightarrow S/xN\longrightarrow V
       \longrightarrow0,\\
 0\longrightarrow V\xrightarrow{\cdot x}S/xN
       \longrightarrow S/xS\longrightarrow0.
 \end{gathered}
\]
The second map is 
injective because
$x$ is a nonzerodivisor in $S$. 
Its image is $xS/xN$. By Whitney's formula
$
 c(N/xN)c(V)=c(S/xN)=c(V)c(S/xS)=c(V),
$
so that$c(N/xN)=1$. Multiplication by
$\varpi^{r-\dd}$ identifies $E|_{Z_r}$ with $N/xN$ as ordinary
bundles. Since the $Z_r$ exhaust $\mathcal G$, this proves $c(E)=1$.

To compute the divided equivariant Chern classes, let
$T=\mathbb G_m$ act by scaling $\varpi$, and let
$\eps\in H_T^2(\mathrm{pt})$ be its standard equivariant parameter.
Since $x=\varpi^n$ has weight $n$, the equivariant version of the
second sequence is
\[
 0\longrightarrow V\otimes\C_{n\eps}\xrightarrow{\cdot x}S/xN
       \longrightarrow S/xS\longrightarrow0.
\]
Here $\C_{n\eps}$ is the constant line bundle with $T$-weight $n$,
whose equivariant first Chern class is $n\eps$. 
Forgetting equivariance makes this line bundle trivial and recovers the ordinary multiplication map above.
If $v_1,\ldots,v_r$ are the equivariant Chern roots of $V$, the equivariant version of Whitney's
formula gives
$
 c_T(N/xN;s)=c_T(S/xS;s)
       \prod_{a=1}^r\frac{1+(v_a+n\eps)s}{1+v_as}, \ c_T(-;s)=\sum_j c_j^T(-)s^j.
$
The weights of $S/xS$ are $0,\eps,\ldots,(n-1)\eps$, so its term
linear in $\eps$ has Chern degree one. Comparing the terms linear
in $\eps$ in higher Chern degrees yields
\begin{equation}\label{eq:bott-chern}
 \left.\frac{c_{i+1}^T(E)}{n\eps}\right|_{\eps=0,Z_r}
   =(-1)^ip_i(V)=(-1)^iu_i^\vee|_{Z_r},\qquad 1\le i<n.
\end{equation}
The division is legitimate on $\mathcal G$: its equivariant cohomology
is torsion-free over $\C[\eps]$, and the positive ordinary Chern
classes of $E$ vanish. 
The translation in \eqref{eq:gap-quot} twists
$E$ by a constant character; it does not change the displayed specializations, again because the positive ordinary Chern classes vanish.

The classes on the left-hand side of \eqref{eq:bott-chern} restrict to the classes
$e_{i+1}$ on $\J$ from \cite{OY17} and hence $u_i^\vee|_\J=(-1)^{i}e_{i+1}$. By \cite[Theorems 1.1.2 and 1.2.1]{OY17},
$H^*(\J)$ is generated by these classes, and $H^*(\tJ)$ by their
pullbacks together with $\xi_1,\ldots,\xi_n$. Thus restriction from
the ambient Grassmannian or affine flag variety is surjective.

We remark that the Chern filtration used in \cite{OY17} is defined first in
localized equivariant full-flag cohomology by polynomial degree in
the flag Chern classes, with coefficients in $\C[\eps,\eps^{-1}]$ \cite[Definitions 8.1.1, 8.2.1]{OY16}.
One then intersects with equivariant cohomology and specializes to
ordinary cohomology. It is this saturated filtration that agrees
with $\mathsf{P}$, by \cite[Proposition 8.5.4]{OY16} and
\cite[Proposition 5.3.1]{OY17}. Since $c_{i+1}^T(E)$ is a polynomial
of degree $i+1$ in the equivariant flag Chern classes, this gives
\begin{equation}\label{eq:chern-bounds}
 u_i^\vee|_{\J}\in \mathsf{P}_{\le i+1}H^{2i}(\J),\qquad
 \xi_j|_{\tJ}\in \mathsf{P}_{\le1}H^2(\tJ).
\end{equation}

\begin{lem}\label{lem:restriction}
For $Y=\J$ or $\tJ$ and $k,r\ge0$, restriction induces a surjection
\[
 \mathsf{P}_{\le k+r}H^{2k}(Y)\twoheadrightarrow H^{2k}(Y_{\le r}).
\]
\end{lem}
\begin{proof}
The gap loci and their complements have compatible affine pavings,
by \cite{Hikita} and \cite[\S8.3, Theorem 8.3]{KTpublished}. In particular,
$H_c^{\mathrm{odd}}(Y\setminus Y_{\le r})=0$. The compactly supported
cohomology sequence therefore gives
$H^{2k}(Y)\twoheadrightarrow H^{2k}(Y_{\le r})$.

We explain the spherical case; the case with flags is similar. Consider the commutative square of restrictions
\[
 \begin{array}{ccc}
 H^*(\mathcal G)&\longrightarrow&H^*(\J)\\
 \downarrow&&\downarrow\\
 H^*(Z_r)&\longrightarrow&H^*(\J_{\le r}).
 \end{array}
\]
The top map is surjective by Chern generation and the right map by
the preceding paragraph. By commutativity also
the bottom map is surjective. Lemma~\ref{lem:bott} and
\eqref{eq:chern-bounds} give the required lifts of Chern/perverse degree
at most $k+r$.
\end{proof}

\begin{proof}[Proof of Theorem \ref{thm:coprime}]
Lemma~\ref{lem:bott} shows that every class in
$H^{2k}(\J_{\le r})$ is the restriction of a linear combination
of ambient cup monomials
\[
 (u^\vee)^\alpha=\prod_{i=1}^{n-1}(u_i^\vee)^{\alpha_i},
 \qquad
 \sum_i i\alpha_i=k,\qquad
 |\alpha|:=\sum_i\alpha_i\le r.
\]
On the other hand, the tautological-bundle computation identifies
the restriction of $u_i^\vee$ to $\J$ with $(-1)^ie_{i+1}$.
Its Chern-filtration degree is therefore at most $i+1$.
By multiplicativity and Oblomkov--Yun's identification of the
Chern and perverse filtrations, the restriction of this monomial
lies in
\[
 \mathsf{P}_{\le\sum_i(i+1)\alpha_i}H^{2k}(\J)
 =
 \mathsf{P}_{\le k+|\alpha|}H^{2k}(\J)
 \subseteq \mathsf{P}_{\le k+r}H^{2k}(\J).
\]
Consequently, restriction induces the surjection of
Lemma~\ref{lem:restriction}:
\begin{equation}\label{eq:sketch-restriction}
 \mathsf{P}_{\le k+r}H^{2k}(\J)
 \twoheadrightarrow H^{2k}(\J_{\le r}).
\end{equation}

The character computations of Oblomkov--Yun
\cite{OY16,OY17} and Berest--Etingof--Ginzburg \cite{BEG},
Hikita's formula transported through the involution of
\cite[Theorem~8.3]{KTpublished}, and Mellit's rational shuffle
theorem with its reciprocal specialization
\cite{Hikita,Mellit,GLWX} give
\begin{equation}\label{eq:sketch-characters}
 \sum_{k,a}z^a\dim\operatorname{gr}^{\mathsf{P}}_aH^{2k}(\J)
 =
 \sum_{k,s}z^{k+s}\dim\operatorname{gr}^{\mathsf{Q}}_sH^{2k}(\J).
\end{equation}
In the classical case $m=n+1$, the singly graded comparison between
the Cherednik module and diagonal coinvariants was established by
Gordon \cite[\S5]{Gordon}, with the grading shifted by $\dd$ to match
our normalization. In the unflagged case, the common value in
\eqref{eq:sketch-characters} is up to a normalization the rational Catalan $z$-number
\[
z^{\delta}C_{n,m}(z,z^{-1})=\prod_{i=2}^{n}\frac{1-z^{m+i-1}}{1-z^i}.
\]
Only the equality in \eqref{eq:sketch-characters} is needed below.

Fix an integer $p$. The surjections
\eqref{eq:sketch-restriction}, with $r=p-k$, imply
$
 \dim \mathsf{P}_{\le p}H^{2k}(\J)
 \ge \dim H^{2k}(\J_{\le p-k}).
$
Taking cumulative coefficients through
$z^p$ in \eqref{eq:sketch-characters} gives
\[
 \sum_k\dim \mathsf{P}_{\le p}H^{2k}(\J)=\sum_{k+s\le p}
       \dim\operatorname{gr}^{\mathsf{Q}}_sH^{2k}(\J)=\sum_k\dim H^{2k}(\J_{\le p-k}).
\]
Every individual difference is nonnegative and their sum is zero.
Hence each difference vanishes, and the surjections
\eqref{eq:sketch-restriction} are isomorphisms.

Finally, by definition,
$
 \mathsf{Q}^{>r}H^{2k}(\J)
 =
 \ker\!\left(
 H^{2k}(\J)\longrightarrow H^{2k}(\J_{\le r})
 \right).
$
Every class in $H^{2k}(\J)$ therefore has a unique representative
in $\mathsf{P}_{\le k+r}H^{2k}(\J)$ with the same restriction to
$\J_{\le r}$. Subtracting that representative gives the direct-sum
decomposition
\[
 H^{2k}(\J)
 =
 \mathsf{P}_{\le k+r}H^{2k}(\J)
 \oplus \mathsf{Q}^{>r}H^{2k}(\J).
\]
For $r<0$, the same dimension comparison gives
$\mathsf{P}_{\le k+r}H^{2k}(\J)=0$, so the assertion remains valid.
Odd cohomology vanishes by the affine paving.

For full flags, replace $Z_r$ by the corresponding flag bundle.
A finite-flag monomial of degree $2a$ contributes at most $a$
to the perverse degree, while the remaining monomial contributes
at most $k-a+r$; thus the bound remains $k+r$.
The same dimension argument applies, with common polynomial the
rational parking function $z$-number
$[m]_z^{n-1}$.

Both filtrations and the restriction maps are compatible with
the Springer action. The result for arbitrary $\J^\nu$ therefore follows by taking
the appropriate Young-subgroup invariants, which is exact over
$\mathbb C$.
\end{proof}
\begin{cor}
The parabolic Hilb-vs-Quot conjecture \cite[Conjecture 2]{KTpublished} and the ORS conjecture \cite[Conjecture 2]{ORS} hold for the singularities $\{y^n=x^m\}\subset \bb{A}_\bb{C}^2$, $(m,n)=1$.
\end{cor}
\begin{proof}
The relation between the $\mathsf{Q}$-filtration and the Quot schemes, as well as the relation between the $\mathsf{P}$-filtration and the Hilbert schemes, is explained in \cite{KTpublished}. A priori, the latter is only explained for the unflagged case. For the parabolic Hilbert schemes, one can use for example \cite[Proof of Theorem 9.7]{GORS}\footnote{Note that the proof there is incomplete. As far as we are aware, there is no published proof of the parabolic Macdonald formula in the literature.}. 
\end{proof}
\begin{ex}[The spherical $y^3-x^5$ case]
Here $R=\C[[\varpi^3,\varpi^5]]$. 
Put $a=u_1^\vee$ and $b=u_2^\vee$, and use
the same letters for their restrictions. Then
\[
 H^*(\mathcal G)=\C[a,b],
 \qquad \deg a=2,\quad\deg b=4.
\]

We recall the affine paving on $\J$. Each valuation semimodule can be written as
$\Delta_b=\bigcup_{i=0}^2(b_i+3\mathbb Z_{\ge0})$, with
$b_i\equiv i\pmod3$. 
The corresponding cell
$A_{\Delta_b}$ parametrizes lattices with valuation set $\Delta_b$. The seven cells are tabulated below.
\begin{center}\small
\begin{tabular}{@{}c c c l@{}}\toprule
$(b_0,b_1,b_2)$ & $\mathbb Z_{\ge0}\setminus\Delta_b$ & $r$ & $\dim A_{\Delta_b}$ \\\midrule
$(6,4,5)$ & $\{0,1,2,3\}$ & $0$ & $0$  \\
$(3,7,5)$ & $\{0,1,2,4\}$ & $1$ & $1$  \\
$(3,4,8)$ & $\{0,1,2,5\}$ & $1$ & $2$  \\
$(6,7,2)$ & $\{0,1,3,4\}$ & $2$ & $2$  \\
$(9,4,2)$ & $\{0,1,3,6\}$ & $2$ & $3$  \\
$(6,1,8)$ & $\{0,2,3,5\}$ & $3$ & $3$  \\
$(0,10,5)$ & $\{1,2,4,7\}$ & $4$ & $4$ \\\bottomrule
\end{tabular}
\end{center}
The table and Theorem~\ref{thm:coprime} give that the bigraded Poincar\'e
polynomials for $\sf Q$ and $\sf P$ are
\begin{align}
 \mathcal Q_\J(q,v):=\sum_{k,r}q^rv^k\dim\gr^{\sf Q}_rH^{2k}(\J)=1+qv+(q+q^2)v^2+(q^2+q^3)v^3+q^4v^4,\\
 \mathcal P_\J(q,v):=\sum_{k,p}q^pv^k\dim\gr^{\sf P}_pH^{2k}(\J)=1+q^2v+(q^3+q^4)v^2+(q^5+q^6)v^3+q^8v^4.
\end{align}
Equivalently, $\mathcal Q_\J(q,v)=v^4C_{5,3}(q,v^{-1})$, where
$C_{5,3}(q,t)=q^4+q^3t+q^2t^2+qt^3+t^4+q^2t+qt^2$ is the
rational $q,t$-Catalan polynomial in this case. 

We can also compute the perverse polynomial independently of
Theorem~\ref{thm:coprime}. 
Define equivariant lifts of $a, b$ as follows:
\[
 \alpha=p_1^T(S/M)-6\eps,\qquad
 \beta=p_2^T(S/M)-7\eps p_1^T(S/M)+28\eps^2.
\]
Restriction to the fixed points is injective in equivariant cohomology,
and the values at fixed points after setting $\eps=1$ are
$(0,0),(1,0),(2,2),(2,-2),(4,4),(4,-4),(8,0)$.
These values verify the
 homogeneous identities
\begin{equation}
\label{eq:equivrels}
 (\alpha-2\eps)(\alpha-4\eps)\beta=0,\qquad
 3\beta^2+\eps\alpha(\alpha-\eps)(\alpha-8\eps)=0.
\end{equation}
Specializing gives $a^2b=b^2=0$, while $a^5=0$ by dimension.
Chern generation and dimension counting now prove
$
 H^*(\J)\cong \C[a,b]/(a^5,a^2b,b^2).
$

By the tautological-bundle computation,
$\alpha,\beta$ are obtained from OY's divided Chern generators by an invertible
triangular change preserving the Chern filtration. We may therefore
compute the localized Chern filtration using the weights
$
 \operatorname{wt}(\alpha)=2,
 \operatorname{wt}(\beta)=3,
 \operatorname{wt}(\eps)=0.
$

Multiplying the first relation in Eq.~\eqref{eq:equivrels} by $\beta$, substituting the second, and
cancelling $\eps$ gives
$
 \alpha(\alpha-\eps)(\alpha-2\eps)
        (\alpha-4\eps)(\alpha-8\eps)=0.
$
Here cancellation is valid because $H_T^*(\J)$ is
$\eps$-torsion-free.
These relations reduce every monomial to a $\C[\eps]$-linear
combination of
$
 1,\ \alpha,\ \beta,\ \alpha^2,\ \alpha\beta,\ \alpha^3,\ \alpha^4,
$
without increasing its Chern weight. Indeed,
$
 \alpha^2\beta=6\eps\alpha\beta-8\eps^2\beta,\
 \beta^2=-\frac{\eps}{3}\alpha^3
          +3\eps^2\alpha^2-\frac83\eps^3\alpha,
$
and the remaining relation expresses $\alpha^5$ in terms of lower
powers of $\alpha$. By equivariant formality,
these classes form a $\C[\eps]$- basis of $H_T^*(\J)$.
Consequently, each weighted filtration step is the
$\C[\eps]$-span of the basis elements of weight at most that step.
It is a direct summand, hence already $\eps$-saturated.
After specializing $\eps=0$,
\cite[Definition 5.2.1 and Proposition 5.3.1]{OY17} show that
$
 1,\ a,\ b,\ a^2,\ ab,\ a^3,\ a^4$
form a basis adapted to the perverse filtration, 
and 
$
 \mathcal P_\J(q,v)
 =1+q^2v+(q^3+q^4)v^2
       +(q^5+q^6)v^3+q^8v^4.
$

Finally, we discuss the role of Lemma \ref{lem:restriction} here. For $r=1$ we have $\J_{\le1}=Z_1\cong\mathbb P^2$. 
If
$h=c_1(\mathcal O(1))$, then $a\mapsto h$ and $b\mapsto h^2$ under restriction.
For $r=2$, the quotient bundle $V$ on $Z_2$ has rank two.
The calculation in \eqref{eq:bott-chern} gives $\left.\frac{c_{i+1}^T(E)}{3\eps}\right|_{\eps=0}=(-1)^ip_i(V)$ for all $i\geq 1$. Since $E$ has rank $3$, $p_i(V)=0$ for $i\geq 3$.
Newton's identities give
$
 2p_3(V)=3ab-a^3,\qquad
 2p_4(V)=b^2+2a^2b-a^4.$
We have the following
commutative square as in Lemma~\ref{lem:restriction}:
$$
\begin{tikzcd}
 \C[a,b]\arrow[r]\arrow[d]&\C[a,b]/(a^5,a^2b,b^2)\arrow[d]\\
 \C[a,b]/(a^3-3ab,b^2-a^2b)\arrow[r]& \C[a,b]/(a^3-3ab,a^2b,b^2).
 \end{tikzcd}
$$
The bottom-left ring has Hilbert series $1+v+2v^2+v^3+v^4$; the bottom-right ring has
$1+v+2v^2+v^3$, agreeing with the five cells of $\J_{\le2}$. 

The perverse filtration computation above and the restriction maps consequently give
\[
\setlength{\arraycolsep}{4pt}
\begin{array}{c|l|c|l}
p & \mathsf{P}_{\le p}H^*(\J)
& r & \mathsf{Q}^{>r}H^*(\J)\\ \hline
p<0 & 0
& r<0 & H^*(\J)\\
0,1 & \langle1\rangle
& 0 & \langle a,b,a^2,ab,a^3,a^4\rangle\\
2 & \langle1,a\rangle
& 1 & \langle a^2-b,ab,a^3,a^4\rangle\\
3 & \langle1,a,b\rangle
& 2 & \langle a^3-3ab,a^4\rangle\\
4 & \langle1,a,b,a^2\rangle
& 3 & \langle a^4\rangle\\
5 & \langle1,a,b,a^2,ab\rangle
& r\ge4 & 0\\
6,7 & \langle1,a,b,a^2,ab,a^3\rangle
& & \\
p\ge8 & H^*(\J)
& &
\end{array}
\]
Here the \(r=1,2\) entries for $\sf{Q}$ follow from the restriction maps computed above. For \(r=3\), the paving shows that restriction is an isomorphism in degrees below \(8\), while \(H^8(\mathcal J_{\le3})=0\); hence its kernel is precisely \(\mathbb C a^4\).

Using the above, it is easy to verify the decomposition
\[
 H^{2k}(\J)
 =\mathsf{P}_{\le k+r}H^{2k}(\J)
  \oplus\mathsf{Q}^{>r}H^{2k}(\J).
\]
For example in degrees $4$ and $6$, we have
\[
\begin{array}{l|l|l|l}
r & \mathsf{P}_{\leq k+r}H^{2k}(\J) & \mathsf{Q}^{>r}H^{2k}(\J) & \mathsf{P}_{\leq k+r}H^{2k}(\J)\oplus \mathsf{Q}^{>r}H^{2k}(\J)\\
\hline
\leq 0 &\mathsf{P}_{\le 2+r}H^4(\J)=0
& \mathsf{Q}^{>r}H^4(\J)=\langle b,a^2\rangle
& H^4(\J)\\
=1 & \mathsf{P}_{\le3}H^4(\J)=\langle b\rangle
& \mathsf{Q}^{>1}H^4(\J)=\langle a^2-b\rangle
& H^4(\J)\\
\geq 2 & \mathsf{P}_{\le 2+r}H^4(\J)=\langle b,a^2\rangle
& \mathsf{Q}^{>r}H^4(\J)=0
& H^4(\J)\\ \hline
\leq 1 & \mathsf{P}_{\le 3+r}H^6(\J)=0
& \mathsf{Q}^{>r}H^6(\J)=\langle ab,a^3\rangle
& H^6(\J)\\
=2 & \mathsf{P}_{\le5}H^6(\J)=\langle ab\rangle
& \mathsf{Q}^{>2}H^6(\J)=\langle a^3-3ab\rangle
& H^6(\J)
\end{array}
\]

\end{ex}

\end{document}